\documentstyle{article}
 \def\vt{t\kern-0.22em\raise.18ex\hbox{\char'47}\lower.18ex\hbox{}\kern-0.08em}
 
\newtheorem{th}{Theorem}[section]
\newtheorem{co}{Corollary}[section]

\newtheorem{ob}{Observation}[section]
\newtheorem{con}{Conjecture}[section]

\newcommand{\old}[1]{{}} 
\newcounter{obr}
\newcounter{tabul}
\begin{document}
\title{Largest reduced neighborhood clique cover number revisited   
}
\author{Farhad Shahrokhi\\
Department of Computer Science and Engineering,  UNT\\
farhad@cs.unt.edu
}

\date{}
\maketitle
\thispagestyle{empty}
\date{} \maketitle
%%%%%%%%%%%%%%%%%%%%%%%%%%%%%%%%%%%%%%%%%%%%%%%%%%%%%%%%%%%%%%%%%%%%%%%%%%%%%%%%%%%%%%%%%%%%%%%%%%%
\begin{abstract} 

 Let  $G$ be a graph and $t\ge 0$.  
 The largest reduced neighborhood clique cover number  of $G$,   denoted by 
${\hat\beta}_t(G)$, 
 is the largest, overall $t$-shallow minors $H$ of $G$, of the smallest number of  cliques that can  cover   any  closed neighborhood of a vertex in $H$.    
 It  is known that  ${\hat\beta}_t(G)\le  s_t$, where $G$ is an incomparability graph and $s_t$ is the  number of leaves in  a largest $t-$shallow minor which is isomorphic to an induced star on $s_t$ leaves.
 In this paper  we give an overview of  the   
  properties of  ${\hat\beta}_t(G)$ including the connections to  the greatest reduced average 
density of $G$, or $\bigtriangledown_t(G)$,  introduce  the class of graphs with bounded neighborhood clique cover number, and derive a simple lower  and an upper bound 
for this important graph parameter.  We announce two  conjectures, one  for the value of ${\hat\beta}_t(G)$, and another for a  separator theorem (with respect to a certain  measure)  for  an interesting  class of graphs, namely the class of  incomparability graphs which we suspect  to have  a polynomial bounded neighborhood clique cover number, when the size of a largest induced star is bounded. 
 
 \end{abstract}

\section{Introduction}
 This paper is a sequel to our paper \cite{Sh5}.  
We assume the reader is familiar with  standard graph theory. Throughout this paper $G=(V,E)$ denotes an undirected graph. Recall that a graph $G$ is $k-$degenerate ($k\ge 0$), if every induced subgraph of $G$ has a vertex of degree at most $k$. Degeneracy of $G$ is the smallest  integer $k$ so that $G$ is $k-$degenerate. Graphs with small degeneracy have  nice structural and algorithm properties. Nesetril and Ossona de Mendez introduced  an important graph parameter  which 
is a generalization of degeneracy. In simple words they introduced the notion of the maximum edge density of a graph  taken overall {\it $t-$shallow minors}. 

A  {\it $t-$shallow minor}, or a {\it $t-$minor} of $G$ in short,  is a minor   of $G$ which is obtained by contracting  connected subgraphs of radius  at most $t$, and deleting vertices (but not edges). Nesetril and Ossona de Mendez  introduced the greatest reduced average density   of $G$ (grad of $G$ in short), or  $\bigtriangledown_t(G)$, to be the maximum edge density of any $t-$minor in $G$. It is easily seen that ${\hat {\delta}(G)\over 2} \le \bigtriangledown_t(G)$, where 
$\hat \delta(G)$ is the degeneracy of $G$.  They  define $G$ to have bounded expansion, 
 if $\bigtriangledown_t(G)$  
is finite for every $t\ge 0$. They  explored very nice structural and algorithmic properties  of the class of  bounded expansion graphs 
that  contains many  traditionally known ``sparse'' graphs, including the class of $H-$minor free graphs \cite{N1,N2,N3,N4}. 

 We introduced the largest reduced neighborhood clique cover number of $G$,  in \cite{Sh5}.  Informally, consider  the minimum number of disjoint cliques that covers the  closed neighborhood of any vertex in a graph;  Now take the maximum value of such a minimum overall $t-$minors of the graph. 
 Formally,  for  a graph $H$,   let $\beta(H)$ denote the clique cover number of $H$, that is, the minimum number of disjoint cliques that partitions $V(H)$.   
 Now  for any $x\in V(H)$,  let $H_x$ denote the 
 the closed neighborhood of $x$ in $H$, and note that $\beta(H_x)\le deg(x)$, where $deg_{H}(x)$ is the degree of $x$ in $H$. Now, 
let ${\tilde \beta(H)}=min_{x\in V(H)}\{\beta(H_X)\}$, and note 
that ${\tilde \beta(H)}\le {\hat \delta(H)}$. 
 Next  for any graph $G$ and $t\ge 0$  define  
 largest reduced neighborhood clique cover number  of $G$, denoted by  ${\hat\beta}_t(G)$ to be  the largest value of $\tilde\beta(H)$ for any $t-$minor $H$ of $G$. 
  We say $G$ has a bounded  neighborhood clique cover number 
  if ${\hat\beta}_t(G)$ has a finite value for each $t\ge 0$.  
   Note that $\hat\beta_t({K_n})={1}$  for any $t\ge 0$, nonetheless 
 $\bigtriangledown_t(K_n)={n-1\over 2}$. Furthermore, one can construct non trivial classes of graphs so that for  every $G$ in the class ${\hat\beta}_t(G)$ is small, that is  bounded by a constant, whereas,  $\bigtriangledown_t(G)$ is arbitrary large.
 For instance, let $G=(V,E)$ be a connected graph which  is the  complement of a bipartite graph, where each partite class has $n$ vertices. Then   ${\hat\beta}_t(G)\le 2$, whereas, $\bigtriangledown_t(G)={|E|\over |V|}\ge {n-1\over 2}$, for any $t\ge 0$. Additionally, for any chordal graph $G$, ${\hat\beta}_t(G)=1$ \cite{Sh5}, but of course one can construct very dense non trivial chordal  graphs   $G$ for which ${\hat\beta}_t(G)$ is unbounded. 
 
 ${\hat\beta}_t(G)$,  is an effective tool  to  study   the properties  of those  graphs that  are not   `` sufficiently sparse'', to  have  a bounded expansion, but yet there is need to explore their properties. 
 For instance, another interesting class of  graphs  for  which  $\hat\beta(G)$ is bounded, but grad of $G$ can be arbitrary large is the intersection graph of fact objects (spheres, cubes, boxes with bounded aspect ratio)\cite{Chan}  when  geometric dimension  is bounded.  Specifically, see  \cite{Sh5} for the following Theorem.
  
\begin{th}\label{t3}
{\sl 
 Let $G$ be the intersection graph of fat objects in $R^d$  (spheres, cubes, boxes with bounded aspect ratio), then, ${\hat\beta}_t(G)=O(b^d.t^{2d})$, where $b$ is a constant that depends on the shape of the object.
 }
 \end{th}
 
  Section two contains    a simple lower  bound and   an upper bound on 
${\hat\beta}_t(G)$ in terms of the clique cover width of $G$,  and some  constructions that measures the ratio of  the upper bound to the  lower bound. Section three  contains two conjectures related to incomparability  graphs that arise from our studies here.

 \section{Bounds on ${\hat\beta}_t(G)$}
 
It is interesting to observe that  ${\hat\beta}_0(K_{n,n})=n$, therefore,  ${\hat\beta}(K_{n,n})$ is not bounded. In fact, the following observation is easy to prove.
   
\begin{ob}\label{o1}
   {\sl  Let $p$ be the largest integer so that a  $t-$shallow minor of $G$ is isomorphic to $K_{p,p}$, then ${\hat\beta}_t(G)\ge p$.  
    }
   \end{ob}

%In this paper we focus on the connection between  ${\hat\beta}_t(G)$, $CCW(G)$, or the %clique cover width of $G$,  and hence the  separation  properties of $G$. 
 
%consequence of the upper bounds derived in 
%section two is that for any incomparability graph $G$, ${\hat\beta}_t(G)=O(t.s)$, %where $s$ %is the maximum number of leaves in an induced star of $G$. In fact , we %show a stronger %result that involves the  concept of the clique cover width 
%of $G$, or $CCW(G)$.

For a clique cover $C$ in $G$, 
the {\it clique cover graph} of $C$ 
is the  graph obtained by contracting the vertices of  each
clique in $C$ into a single vertex. The {\it clique cover width} of  $G$,
denoted by $CCW(G)$, is the minimum  value of the  bandwidth of 
all clique cover graphs in $G$\cite{Sh1,Sh2,Sh3}.   
In this paper when we  write  $C=\{C_1,C_2,....,C_k\}$,
we mean $C$ is an ordered set. Let $ab$ be an edge width $a\in C_i$ and $b\in C_j, j>i$, and  let 
$W(e)=j-i$. We call $W(e)$ the {\it width} of $e$. 
An important  application of the clique cover width 
is in the derivation of separation theorems in dense graphs, where separation can be
defined for  other types {\it measures }\cite{Sh1}, instead of just the number of vertices.
Recall that according to  the planar separation  theorem, any $n$ vertex planar graph can be separated into  two subgraphs, each having at most $2n/3$ vertices,   by removing $O({\sqrt n})$ vertices.
Any $G$  can be separated  with respect to an optimal (or feasible) 
set $C$of cliques (utilizing $CCW(G)$): There is partition of $\{A,S,B\}$ of $V(G)$ so that $(i)$ there are no edges between $A$ and $B$, $(ii)$ $S$ can covered with at most $CCW(G)$ many cliques from $C$, and $(iii)$ $A$ and $B$ are each covered with at most $2|C|/3$ cliques from $C$\cite{Sh1,Sh3}. 
 
%An important  tool for the design of a divide and
%conquer algorithm is separation.

%\begin{ob}\label{o2}
 %{\sl   For any graph $H$
%${\tilde \beta(H)}\le CCW(H)$, and hence ${\hat\beta}_t(G)\le k$, where $k$  is  the largest %clique cover width   of any $t-$shallow minor of $G$. 
 %}
%\end{ob}

%Let $\delta(G)$ denote the minimum degree of $G$. Recall that  the degeneracy of %$G$, denoted by ${\hat\delta}(G)$ is the largest minimum degree among all %induced subgraphs of $G$. 

 %Note  that ${\hat \delta(G)\over 2}\le \bigtriangledown_0(G)\le\bigtriangledown_t(G)$. 

 %\section{Main Results}

   \begin{th}\label{t1}
   {\sl  
   For any graph $G$, ${\hat\beta}_t(G)\le k+1$, where $k$  is  the largest clique cover width   of any $t-$shallow minor of $G$. 
   
    }
    \end{th}
    {\bf Proof.} Let $\{C_1,C_2,....,C_K\}$ be a clique cover of  a graph $H$.  Let $e_a=ab, a\in C_1, b\in C_i$ be an edge of largest width incident to $a$.  
   Let $e^*$ be an edge having an end point in $C_1$
    with  $W(e^*)=\min\{W(e_a)|a\in C_1\}$. 
   By definition of $e^*$, $H_a$ can be covered with $W(e^)*+1$ cliques, and hence ${\tilde \beta(H)}\le W(e^*)+1$.  Therefore ${\tilde \beta(H)}\le CCW(H)+1$, since $CCW(H)\ge W(e^*)$. 
   To finish the proof take $H$ to be a $t$-minor of $G$.  
    $\Box$.
    
    \begin{co}\label{c1}
    {\sl Let $k$  denote   the largest clique cover width   of any $t-$shallow minor of $G$, and $p$ be largest integer so that  any $t-$shallow minor of $G$ is isomorphic to $K_{p,p}$. 
    Then, $p\le {\hat\beta}_t(G)\le k+1$.
    
    }
   \end{co}
    
It is easy  to verify that $CCW(H)\le CCW(G)$, for any induced subgraph $H$ of $G$. 
Nonetheless, for a $t-$minor $H$ of $G$, $CCW(H)$, or $k$ in Corollary,  \ref{c1}  may be much larger than $CCW(G)$.  
  Generally speaking, it would nice to know how large the  ratio $k/p$ may be. 
 
 \begin{ob}\label{o2}
    {\sl For any $t\ge 0$, and $n>t$,  there is an $n$ vertex  graph $G$, with $CCW(G)=1$, so that for a $t-$minor 
    $H$ of $G$, $t\ge CCW(H)\ge t/2$. Moreover, in this case, neither $G$, nor $H$ contain 
    $K_{2,2}$ as an induced subgraph.  
    
    }
  \end{ob} 
  
{\bf Justification.} Let $P_n$ be a path of $n$ vertices on vertex set $X=\{x_1,x_2,...,x_n\}$. Now let $S$ be a an independent set of 
   $n$ vertices. To construct $G=(V,E)$ place a perfect matching of cardinality $n$ between 
   $S$ and $X$. It is easily verified that $CCW(G)=1$. Now for a given $n\ge t\ge 0$, contract 
   $x_1,x_2,...,x_t$ into one single vertex to obtain a $t-$minor $H$.  
Observe that $H$ has an induced star on $t$ vertices. Thus, $CCW(H)\ge t/2$ \cite{Sh2}.  Furthermore, it is not difficult to see that $G$ is an incomparability graph (a graph whose complement has a transitive orientation on edges),  and so is $H$, since  $H$ is obtained by  contractions of edges in $G$. Since $H$ is an incomparability graph we must have  $CCW(H)\le s$, where  $s$ is the number of leaves in a largest induced star \cite{Sh2}.
 Finally,
it is easy to verify that neither $H$ or $G$ have $K_{2,2}$ as a subgraph, since $G$ is acyclic. $\Box$

 %Likewise  the lower bound of  $p$ in the above result  can be large  compared to the size a %largest induced complete bipartite subgraph of $G$.   
 
 \begin{ob}\label{o3}
  {\sl  For any $t\ge 0$, and $n>>t$, there is a graph $G$, on $n+t(t+1)$ vertices  that excludes  $K_{2,2}$ as an induced  subgraph, but has a $t-$minor $H$ that contains 
  $K_{t+1,t+1}$ as an induced subgraph. Moreover, $CCW(G)\ge n/2$.}
\end{ob}

{\bf Justification.}  Let $V(G)=A\cup_{i=1}^{t+1}B_i$, where $A$ is a independent set of size $t+1$, and for $i=1,2,...,t$ each $B_i$ is path on $t+1$ vertices; $B_{t+1}$ is a cycle on $n$ vertices.  Now for each $i=1,2,...,t+1$ add a perfect matching of size $t$ between vertices in $A$ and vertices in $B_i$. Thus each vertex in $A$ has degree $t$, where for $i=1,2,...,t$, each vertex of $B_i$ has degree at most $3$. 
   Note that $G$ does not have $K_{2,2}$ as an induced subgraph. Furthermore, since $B_{t+1}$ is a cycle of $n$ vertices, we have  $CCW(G)\ge n/2$. 
   Now for $i=1,2,...,t$, contract each path  $B_i$ into a single vertex. For $B_{i+1}$ contract the first $t+1$ vertices to a vertex.  Then the resulting graph $H$  has an induced subgraph isomorphic to  $K_{t+1,t+1}$. $\Box$

 \section{ Incomparability graphs}
  Recall that a chordal graph does not have any chord-less cycles \cite{Go}.
  An incomparability  graph is a graph whose complement has a transitive orientation \cite{Tr}.  
 Incomparability  graphs are perfect, have geometric realizations, and have recently been subject to intense investigations, due to their intimate connections to string graphs.
 One  wonders if  there is a meaningful converse to Observation $\ref{o1}$. 
 That is,  can one find a suitable  upper bound on ${\hat\beta}_t(G)$ that is related to the lower bound in \ref{o1}? 
 It is less likely that this is the case for all graphs, nonetheless,  we suspect that there is a weak   converse  to \ref{o1} when $G$ is an incomparability graph. Specifically, we have shown that   if an incomparability graph $G$  does not have a $t-$shallow minor  which is isomorphic to an induced star on $s_t$ leaves, then,   ${\hat\beta}_t(G)\le  s_t$. 
 Moreover, we have shown that for any incomparability graph $G$, ${s\over 2}\le CCW(G)\le s$, where $s$ is the number of leaves in a  largest  induced star in $G$. Hence, a natural question  is how large $s_t/s$ can be?

 %Particularly, we prove that for any incomparability graph $G$, and any $t-$shallow minor %$H$ of $G$, we have $CCW(H)=O(t.CCW(G))$; Combining  this result with our result in %\cite{sh1}, we obtain, ${\hat\beta}_t(G)=O(t.s),$ where $s$ is the number of leaves in a %largest induced star in $G$. 

   \begin{con}\label{Con1}
   {\sl Let $G$ be an incomparability graph that does not have an induced star which is isomorphic to an induced star on $s$ leaves. Then, the size of a largest 
   induced star $s_t$ in any $t-$shallow minor of $G$ is at most $O(t.s)$. Consequently, 
    ${\hat\beta}_t(G)=O(t.s)$, for any $t\ge 0$.
   
 }

   \end{con}

 If the  above conjecture were to be true, then  ${\hat\beta}_t(G)=O(t.s)$, where $t$ is the number of leaves in a  largest induced  star in $G$. Note that the conjecture implies that the class of incomparability graphs have a linearly  bounded neighborhood clique cover number, when the size of a largest induced star is fixed. 
 
By observation \ref{o1}, ${\hat\beta}_t(G)\ge p_t$, where  $p_t$ is  the largest integer so that  $K_{p_t,p_t}$ is a $t-$shallow  minor of $G$.  
Hence o get a good  estimate for  ${\hat\beta}_t(G)$ (if the conjecture were to be true), one has to investigate how large $t.s/p_t$  can be.

  It is easy to observe that If $G$ is a chordal  graph, then, ${\hat\beta}_t(G)=1$ \cite{Sh5}. 
  Moreover, the separation property with respect to cliques holds for any chordal graph $G$,  regardless of the value the clique cover width.  Particularly, given a clique tree \cite{Go} of $G$ associated with a set $C$ of maximal cliques, there is one clique $B$ in $C$, so that 
 after removal of $B$, 
each the two remaining (separated) subgraph of $G$ can be covered by at most ${2|C|/3 }$ cliques from $C$. Now let $G$ be an interval  graph; Since  $G$ is  chordal  ${\hat\beta}_t(G)=1$, and additionally $G$ has the stated separation property. Particularly, note that  $G$ is chordal and  also an incomparability graph that does not have a $K_{2,2}$ as an induced subgraph.  In fact,  no $t-$minor of  an interval graph $G$ can have $K_{2,2}$ as an induced subgraph. So one can suspect that if a incomparability  graph $G$ does not have a large $K_{p,p}$ as a $t-$minor, then, $G$ has 
*nice* separation properties with respect to cliques.

 \begin{con}\label{Con2}
   {\sl Let $p$ be fixed, and let $G$ be an incomparability graph that does not have $K_{p,p}$ as a $t-$shallow minor.  Then, there is a clique cover $C$ in $G$ so that the  removal of 
   $O({\sqrt{|C|}})$ cliques from $C$, separates $G$ into two subgraphs so that each subgraph can be covered with at most ${2|C|/3}$ cliques from $C$.

 }   
\end{con}
   
 We remark that by a general result of Fox and Pach \cite{FP} (see also  an  earlier result of Bodlaender and  Thilikos on $k-$chordal graphs \cite{BT}),  any incomparability graph $G$ on $n$ vertices and $m$ edges  has a separation  $(L,S,R)$ so that   
$S=O({\sqrt{m}})$, and $|L|,|S|\le 2n/3$, but conjecture  \ref{Con2} does not follow from these  result.

\end{document}